\documentclass[12pt]{article}
\usepackage{amssymb,latexsym,amsmath}
\usepackage{amsthm}
\usepackage{geometry}
\usepackage{graphicx,color}
\usepackage{lineno}
\theoremstyle{plain}

\theoremstyle{definition}

\numberwithin {equation}{section}
\usepackage{longtable}

\usepackage{float}
\restylefloat{table}

\begin{document}

\title
{A dynamic Stackelberg game for green supply chain management}
\author{Mehrnoosh Khademi $^a$\thanks{mehrnushkhademi@yahoo.com} \and Massimiliano Ferrara $^{b,d}$\thanks{massimiliano.ferrara@unirc.it} \and
Mehdi Salimi $^{c,d}$\thanks{Corresponding author:
mehdi.salimi@tu-dresden.de \& mehdi.salimi@medalics.org}\and
Somayeh Sharifi $^e$\thanks{s.sharifi@iauh.ac.ir}}
\date{}
\maketitle
\begin{center}
$^{a}$Department of Industrial Engineering, Mazandaran University
of Science and Technology, Babol, Iran \\ $^{b}$Department of Law
and Economics, Universit\`{a} Mediterranea di Reggio Calabria,
Italy \\ $^{c}$Center for Dynamics, Department of Mathematics,
Technische Universit{\"a}t Dresden, Germany \\ $^{d}$MEDAlics,
Research Center at Universit\`{a} per Stranieri Dante Alighieri,
Reggio Calabria, Italy \\ $^{e}$Young Researchers and Elite Club,
Hamedan Branch, Islamic Azad University, Hamedan, Iran
\end{center}
\maketitle

%\date{}
\begin{abstract}\noindent
In this paper, we establish a dynamic game to allocate CSR
(Corporate Social Responsibility) to the members of a supply
chain. We propose a model of a three-tier supply chain in a
decentralized state which includes a supplier, a  manufacturer and
a retailer. For analyzing supply chain performance in
decentralized state and the relationships between the members of
the supply chain, we use a Stackelberg game and consider in this
paper a hierarchical equilibrium solution for a two-level game. In
particular, we formulate a model that crosses through
multi-periods with the help of a dynamic discrete Stackelberg
game. We obtain an equilibrium point at which both the profits of
members and the level of CSR taken up by supply chains is
maximized.


 \end{abstract}
%---------------------------------------------------------------------------------------%
\label{sec:Introduction}\section{Introduction} In recent years,
companies and firms have been showing an ongoing interest in favor
of CSR. This is mainly because of increasing consumer awareness of
several CSR issues, e.g.\ the environment, human rights and
safety. In addition, the firms are also forced to accept CSR due
to government policies and regulations. Recently CSR has gained
recognition and importance as field of research field
\cite{Cetindamar,Joyner}. However, the research field still lacks
a consistent definition of CSR and this has been the center of
discussion since several decades. Dahlsrud \cite{Dahlsrud}
presented an overview of different definitions of CSR and
summarized the number of dimensions included in each definition.
There is a positive correlation between CSR and profit
\cite{Orlitzky,Preston}. Moreover, CSR is an effective tool for
supply chain management, for coordination, purchasing,
manufacturing, distribution, and marketing functions
\cite{Hervani}. According to previous studies, the long term
investment on CSR is beneficial for a supply chain. Furthermore, a
sustainable supply chain requires consideration of the social
aspects of the business \cite{Svensson}. Carter et al.
\cite{Carter1} established an effective approach and demonstrated
that environmental purchasing is significantly related to both net
income and cost of goods sold. Carter et al. \cite{Carter} also
pointed towards the importance of CSR in the supply chain, in
particular the role played by the purchasing managers in socially
responsible activities and the effect of these activities on the
supply chain. Sethi \cite{Sethi} introduced a taxonomy in which a
firm's social activities include social obligations as well as
more voluntary social responsibility. And, Carroll
\cite{Carroll1,Carroll} developed a framework for CSR that
consists of economic, legal, and ethical responsibilities.

The members of a supply chain take their decisions based on
maximizing their individual net benefits. In addition, when they
need to accept CSR; this situation leads to an equilibrium status.
Game theory is one of the most effective tools to deal with this
kind of management problems.

A growing number of research papers use game theoretical
applications in supply chain management. Cachon et al.\
\cite{Cachon} discuss Nash equilibrium in noncooperative cases in
a supply chain with one supplier and multiple retailers. Hennet et
al.\ \cite{Hennet} presented a paper to evaluate the efficiency of
different types of contracts between the industrial partners of a
supply chain. They applied game theory methods for decisional
purposes. Tian et al.\ \cite{Tian} presented a system dynamics
model based on evolutionary game theory for green supply chain
management.

In this paper, we consider a discrete time version of the dynamic
differential game. The optimal control theory is the standard tool
for analyzing the differential game theory \ \cite{Kamien}. There
are two different types of information structures in a
differential game, open-loop and feedback information structures.
In an open-loop strategy, the players choose their decisions at
time t, with information of the state at time zero. In contrast,
in a feedback information structure,
 the players use their knowledge of
the current state at time t in order to formulate their decisions
at time $t$ \cite{He}. We formulate a model and study the behavior
for decentralized supply chain networks under CSR conditions with
one leader and two followers. The Stackelberg game model is
recommended and applied here to find an equilibrium point at which
the profit of the members of the supply chain is maximized and the
level of CSR is adopted in the supply chain. We develop an
open-loop Stackelberg game by selecting the supplier as the leader
and both the manufacturer and the retailer as the followers. Using
this approach, the supplier as a leader, can know the optimal
reaction of his followers, and utilizes such processes to maximize
his own profit. The manufacturer and the retailer as followers,
try to maximize their profits by considering all the conditions.
Our model has two levels, at the first level the manufacturer is
the leader and the retailer is the follower and we find the
equilibrium point. At the second level, we consider the supplier
as the leader and the manufacturer as the follower. In fact, we
substitute the response functions of the follower into the
objective function of the leader and we find the final equilibrium
point. We propose a Hamiltonian matrix to solve the optimal
control problem to obtain the equilibrium in this game. The paper
is organized as follows: Section 2 is devoted to the problem
description and assumptions. Objective functions, constraints and
solution of the game are illustrated in Section 3. A numerical
example is shown in Section 4 and we close with a conclusion in
Section 5.

\label{sec:Problem Description and Assumptions}\section{Problem
Description and Assumptions} We consider a three stage Stackelberg
differential game which has three players playing the game over a
fixed finite horizon model. This model is a three-tier,
decentralized vertical control supply chain network (Figure 1).
All retailers and suppliers at the same level make the same
decision. Therefore, consequently the model has only one supplier,
one manufacture and one retailer. The simplified model is shown in
Figure 2.

The dynamic game goes through multi-periods as a repeated game
with complete information. This model has a state variable and
control variables like any dynamic game. We define the state
variable as the level of social responsibility taken up by
companies, and the control variables are the capital amounts
invested while fulfilling the social responsibility. Specifically,
all of the social responsibilities taken up by the firm $j$ at
period $t$ can be expressed as the investment $I_t^j$. We suppose
that $x_t$ evolves according to the following rule:
${x_{t+1}=f(I_t,x_t)}$.
\begin{center}
\begin{figure}[h!]
\begin{minipage}[b]{0.5\linewidth}
\centering
\includegraphics[width=\textwidth]{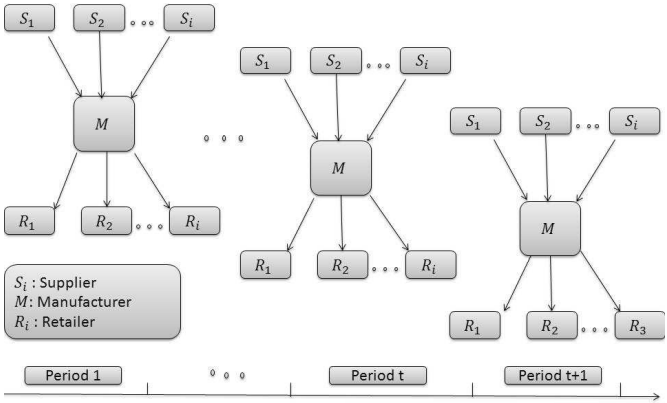}
\caption{Three-tier supply chain network} \label{fig1}
\end{minipage}
\end{figure}
\end{center}

\begin{center}
\begin{figure}[h!]
\begin{minipage}[b]{0.5\linewidth}
\centering
\includegraphics[width=\textwidth]{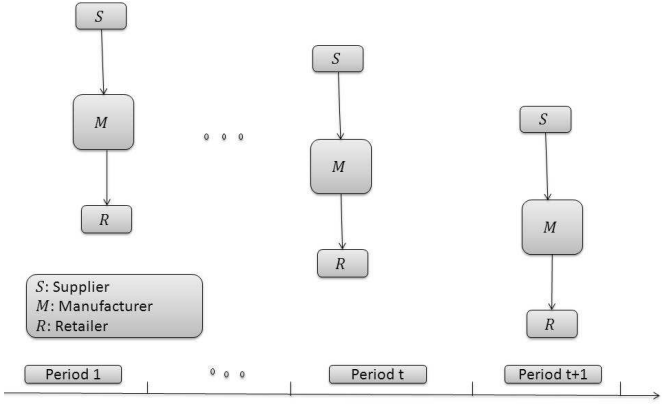}
\caption{Simplified model of three-tier supply chain network}
\label{fig1}
\end{minipage}
\end{figure}
\end{center}

More specifically we have the following assumptions:

The function $B_t(x_t)=\delta x_t$ represents social benefit which
is proportional to social responsibility taken up by the supply
chain system \cite{Batabyal}.

The function $T_t=\tau I_t\big[1+\theta(I_t)\big]$ measures the
value of the tax return to the members of the supply chain
\cite{Feibel}. Both $\tau$ and $\theta$ are tax return policy
parameters. Specifically, $\tau$ is the rate of individual post
tax return on investment (ROI), and $\theta$ is the rate of supply
chain's post tax return on investment (ROI).

The market inverse demand is $P^M(q_t)=a-bq_t$ \cite{Mankiw}.

The accumulation of the level of social responsibility taken up by
the firms is given by $x_{t+1}=\alpha
x_t+\beta_1I_t^S+\beta_2I_t^{M}+\beta_3I_t^{R}$.

Here, ${\beta_1}$ is the rate of converting the supplier's capital
investment in CSR to the amount of CSR taken up by the supply
chain; ${\beta_2}$ is the rate of converting the manufacturer's
capital investment in CSR to the amount of CSR taken up by the
supply chain and ${\beta_3}$ is the rate of converting the
retailer's capital investment in CSR to the amount of CSR taken up
by the supply chain as well \cite{Shi}.

\label{sec:The General Framework}\subsection{The General
Framework} He et al.\ \cite{Guti} illustrate an open-loop
Stackelberg differential game model over a fixed finite horizon
time as detailed in the following:

The follower's optimal control problem is:
\begin{equation}
\begin{split}
&Max_{r(\cdot)}\Big\{J_R(X_0,r(\cdot);w(\cdot))=\int_0^{t}e^{-\rho
t}\pi_R(X(t),w(t),r(t))dt\Big\},\\
\end{split}
\end{equation}
subject to the state equation
\begin{equation}
\begin{split}
&\dot X{(t)}=F(X(t),w(t),r(t)),\\
&X(0)=X_0.
\end{split}
\end{equation}
where the function $F$ represents the rate of sales, $\rho$ is the
followers's discount rate, and $X_0$, is the initial condition.
The follower's Hamiltonian is
\begin{equation}
H_{R}(X,r,\lambda_R,w)=\pi_R(X,w,r)+\lambda_RF(X,w,r),
\end{equation}
where $\lambda_R$ is the vector of the shadow prices associated
with the state variable $X$; and it satisfies the adjoint equation
\begin{equation}
\dot \lambda_R=\rho \lambda_R-\dfrac{\partial
H_R(X,r,\lambda_R,w)}{\partial X}, \quad \lambda_R(T)=0.
\end{equation}
The necessary optimality condition for the follower's problem
satisfies
\begin{equation}
\dfrac{\partial H_R}{\partial r}=0 \Longrightarrow \dfrac{\partial
\pi_R(X,w,r)}{\partial r}+\lambda_R \dfrac{\partial
F(X,w,r)}{\partial r}=0.
\end{equation}
We derive the follower's best response $r^*(X,w,\lambda_R).$

The leader's problem is
\begin{equation}
\begin{split}
&Max_{w(\cdot)}\Big\{J_M(X_0,w(\cdot))=\int_0^{t}e^{-\mu
t}\pi_M(x,w,r(x,w,\lambda_R))dt\Big\},\\
&\dot X=F(X,w,r(x,w,\lambda_R)),\\
&X_0{(0)}=X_0,\\
&\dot\lambda_R=\rho \lambda_R-\dfrac{\partial
H_R(x,r(x,w,\lambda_R),\lambda_R,w)}{\partial x}, \quad
\lambda_R(T)=0,
\end{split}
\end{equation}
where $\mu$ is the leader's discount rate and the above
differential equations are obtained by substituting the follower's
best response $r^*(X,w,\lambda_R)$ into the state equation and the
adjoint equation of the follower, respectively. We formulate
the leader's Hamiltonian as follows
\begin{equation}
H_M=\pi_M(x,\lambda_R,w,r(X,w,\lambda_R),\lambda_M,\varphi)+\lambda
F(X,w,r(X,w,\lambda_R))-\mu\dfrac{\partial
H_R(X,r(X,w,\lambda_R),\lambda_R,w)}{\partial X},
\end{equation}
where $\lambda_M$ and $\mu$ are the shadow associated with $X$ and
$\lambda_R$, respectively, and they satisfy the adjoint equations
\begin{equation}
\begin{split}
\dot \lambda_M&= \mu\lambda_M-\dfrac{\partial
H_M(X,\lambda_R,w,r(X,w,\lambda_R),\lambda_M,w)}{\partial X}\\
&=\mu\lambda_M-\dfrac{\partial
\pi_M(X,w,r(X,w,\lambda_R))}{\partial X}\\
&-\lambda_M\dfrac{\partial F(X,w,r(X,w,\lambda_R))}{\partial
X}-\mu\dfrac{\partial^2
H_R(x,r(X,w,\lambda_R),\lambda_R,w)}{\partial X^2},\\
\dot \varphi&= \mu\varphi-\dfrac{\partial
H_M(X,\lambda_R,w,r(X,w,\lambda_R),\lambda_M,\varphi)}{\partial
\lambda_R}\\
&=\mu\varphi-\lambda_M\dfrac{\partial
F(x,w,r(X,w,\lambda_R))}{\partial \lambda_R}-\mu\dfrac{\partial^2
H_R(X,r(X,w,\lambda_R),\lambda_R,w)}{\partial \lambda \partial X},
\end{split}
\end{equation}
where $\lambda_M(T)=0$ and $\varphi(0)=0$ are the boundary
conditions.\\
We apply the algorithm of the above general model as part of our model.

\label{sec:Notations and Definitions}\subsection{Notations and
Definitions}

 To facilitate the model, certain
parameters and decision variables are used.\\

Table \ref{t1} shows notations and definitions that we use in our
model.
\begin{table}[!ht]
\begin{center}
\begin{tabular}{l l}
Variables & \\
\hline
$t$ & Period $t$ \\
$T$ & Planning horizon\\
$q_t$ & Demand quantity at period $t$\\
$a$ & Market potential\\
$b$ & Price sensitivity\\
$x_t$ & State variable, degree of taking SR\\
$H^S$ & Hamiltonian function of the supplier\\
$H^M$ & Hamiltonian function of the manufacturer\\
$H^R$ & Hamiltonian function of the retailer\\
$J_t^S$ & Objective function of the supplier\\
$J_t^M$ & Objective function of the manufacturer\\
$J_t^R$ & Objective function of the retailer\\
$B^M(x_t)$ & Social benefit of the supplier\\
$B^S(x_t)$ & Social benefit of the manufacturer\\
$B^R(x_t)$ & Social benefit of the retailer\\
$T^S(x_t)$ & Tax return of the supplier\\
$T^M(x_t)$ & Tax return of the manufacturer\\
$T^R(x_t)$ & Tax return of the  retailer\\
$I_t^{M}$ & The amount of investment done by the manufacturer \\
$I_t^S$ & The amount of investment done by the supplier\\
$I_t^R$ & The amount of investment done by the retailer\\
$d$ & The percentage of investment of the supplier payoff\\
$\widehat{d}$ & The percentage of investment of the manufacturer payoff\\
$w$ & Price of the supplier's raw material\\
$z$ & Price of sold product by the retailer\\
$\delta$ & Parameter of the supplier's social benefit\\
$\widehat{\delta}$ & Parameter of the manufacturer's social benefit\\
$\widehat{\widehat{\delta}}$ & Parameter of the retailer's social benefit\\
$\lambda$ & Quantity discount parameter of the price of raw material\\
$\alpha$ & Deteriorating rate of the level of current social responsibility\\
$\tau$ & The rate of individual post tax return on investment (ROI)\\
$\theta$ & The rate of supply chain's post tax return on investment (ROI) \\
$\beta_1$ & The rate of converting the supplier's capital investment in CSR\\
 &to the amount of CSR taken up by the supply chain\\
$\beta_2$ & The rate of converting the manufacturer's capital investment in CSR\\
 &to the amount of CSR taken up by the supply chain\\
$\beta_3$ & The rate of converting the retailer's capital investment in CSR \\
&to the amount of CSR taken up by the supply chain\\

 \hline
\end{tabular}
\end{center}
\vspace*{-3ex} \caption{Notations and Definitions. \label{t1}}
\end{table}
%-------------------------------------------------------------------------------------------------%
\newpage
\label{sec:Objective Functions and Constraints}\section{Objective
Functions and Constraints}

The objective functions are made to depend on the control vectors
and the static variable. The members of the supply chain attempt
to optimize their net profits, which includes minimizing the cost
of raw materials and investment in social responsibility, and
maximizing sale revenues and benefits from taking social
responsibility as well as tax returns. Thus, the objective
function of the supplier is

\begin{equation*}
\begin{split}
J^S&=
\sum_{t=1}^TP_t^Sq_t-cq_t+B_t^S(x_t)+T_t^S(I_t^S,I_t)-I_t^S+dI_t^{M}\\
&=\sum_{t=1}^T wq_t-cq_t+\delta x_t^2+\tau
I_t^S[1+\theta(I_t^S+I_t^{M}+I_t^{R})]-I_t^S+dI_t^{M},
\end{split}
\end{equation*}
where, $P_t^S$ is the price of the supplier's raw material.
${P_t^S}=w$. $B_t^{S}(x_t)$ is the social benefit of the supplier,
$\delta$ is the parameter of the supplier's social benefit and
$T_t^{S}(I_t^{S}, I_t)$ is the tax return of the supplier.
Similarly, the objective function of the manufacturer is
\begin{equation*}
\begin{split}
J^M&=
\sum_{t=1}^TP_t^M(q_t)q_t-P_t^Sq_t+B^M(x_t)+T^M(I_t^M, I_t)-I_t^M+\widehat{d}I_t^R\\
&= \sum_{t=1}^T(a-bq_t)q_t-wq_t+\widehat{\delta}x_t^2+\tau
I_t^M(1+\theta (I_t^S+I_t^M+I_t^R))-I_t^M+\widehat{d}I_t^R,
\end{split}
\end{equation*}
where $P_t^M(q_t)$ is the retail price of the product of the
manufacturer. $B_t^{M}(x_t)$ is the social benefit of the
manufacturer, $\widehat{\delta}$ is the parameter of the
manufacturer's social benefit and $T_t^{M}(I_t^{M}, I_t)$ is the
tax
return of the manufacturer.\\
The objective function of the retailer is
\begin{equation*}
\begin{split}
J^R&=
\sum_{t=1}^TP_t^Rq_t-P_t^M(q_t)q_t+B^R(x_t)+T^R(I_t^R, I_t)-I_t^R\\
&=
\sum_{t=1}^Tzq_t-(a-bq_t)q_t+\widehat{\widehat{\delta}}x_t^2+\tau
I_t^R(1+\theta (I_t^S+I_t^M+I_t^R))-I_t^R,
\end{split}
\end{equation*}
where $P_t^{R}$ is the price at which the retailer sells the
product to the consumer. $P_t^{R}=Z$. $B_t^{R}(x_t)$ is the social
benefit of the retailer, $\widehat{\widehat{\delta}}$ is the
parameter of the retailer's social benefit and $T_t^{R}(I_t^{R},
I_t)$ is the tax return of the retailer.
%----------------------------------------------------------------------------------------%
\label{sec:Mathematical Model: Level One}\subsection{Mathematical
Model: Level One} At this level, we establish a Stackelberg game
between the manufacturer as the leader and the retailer as the
follower. To calculate the equilibrium at this level, first we
calculate the best reaction function of the retailer, then we
determine the manufacturer's optimal decisions based on the
retailer's best reactions.

Since our dynamic differential game is an optimal control problem,
we can apply the Hamiltonian function to find the
equilibrium of the game \cite{Sethi}.\\
Suppose the time interval is $[1,T]$. For any fixed $I_t^S$ and
$I_t^M$ the retailer solves
\begin{equation*}
\begin{split}
&\arg \max_{I_t^R}
\sum_{t=1}^TP_t^Rq_t-P_t^M(q_t)q_t+B^R(x_t)+T^R(I_t^R,
I_t)-I_t^R,\\
\end{split}
\end{equation*}
subject to $x_{t+1}=\alpha x_t+\beta_1I_t^S+\beta_2I_t^{M}+\beta_3I_t^{R}$.\\

We define the retailer's  Hamiltonian for fixed $I_t^S$ and
$I_t^M$ as\begin{equation*}
\begin{split}
H_t^{R}&=J_t^{R}+P_{t+1}^{R}(x_{t+1}).\\
\end{split}
\end{equation*}

By using the conditions for a maximization of this Hamiltonian, we
compute:

\begin{equation}\label{a0}
\begin{split}
I_t^R=\dfrac{1-P_{t+1}^R\beta_3-\tau
\theta(I_t^M+I_t^S)-\tau}{2\theta \tau}.
\end{split}
\end{equation}
The equation of $I_t^R$ which depends on $I_t^S$ and $I_t^M$, says that for any
given strategy of $I_t^S$ and $I_t^M$, there is a unique
optimal response $I_t^R$.
\begin{equation}\label{a1}
x_{t+1}=\dfrac{\partial H_t^{R}}{\partial P_{t+1}^{R}}=\alpha
x_t+\beta_1I_t^S+\beta_2I_t^{M}+\beta_3 I_t^{R},
\end{equation}
and by substituting (\ref{a0}) in (\ref{a1}), we obtain
\begin{equation}\label{a2}
x_{t+1}=(\beta_1-\beta_3/2)I_t^S+(\beta_2-\beta_3/2)I_t^M+\alpha
x_t+\beta_3\dfrac{1-P_{t+1}^R\beta_3-\tau}{2\tau \theta}.
\end{equation}
We also have
\begin{equation}\label{a3}
\begin{split}
P_{t}^R&=\dfrac{\partial H_t^{R}}{\partial
x_{t}}=2\widehat{\widehat{\delta}}x_t+\alpha P_{t+1}^R.
\end{split}
\end{equation}
The above sets of equations define the reaction function of the
retailer.

For any fixed $I_t^S$ the manufacturer solves
\begin{equation*}
\begin{split}
&\arg \max_{I_t^M}
\sum_{t=1}^TP_t^M(q_t)q_t-P_t^Sq_t+B^M(x_t)+T^M(I_t^M, I_t)-I_t^M+\widehat{d}I_t^R,\\
\end{split}
\end{equation*}
subject to $x_{t+1}=\alpha x_t+\beta_1I_t^S+\beta_2I_t^{M}+\beta_3I_t^{R}$.\\

Now, we substitute the value of $I_t^R$ in (\ref{a0}) into $J_t^M$,
and we obtain
\begin{equation}\label{aaa}
\begin{split}
J_t^M&=(a-b q_t)q_t-w q_t+\widehat{\delta}x_t^2+\dfrac{\tau
\theta-\widehat{d}}{2}I_t^S
I_t^M+\dfrac{\tau-1-\beta_3P_{t+1}^R-\widehat{d}}{2}I_t^M-\dfrac{\widehat{d}}{2\theta}.
\end{split}
\end{equation}
The Hamiltonian function of the manufacturer for fixed $I_t^S$ is

\begin{equation}\label{a4}
\begin{split}
H_{t}^M&=J_t^M+P_{t+1}^M(x_{t+1})+u_t(P_t^R),
\end{split}
\end{equation}

consequently, we can obtain the unique optimal response of the
follower from the equations as follows.

\begin{equation*}
\begin{split}
\dfrac{\partial H_t^{M}}{\partial
I_{t}^{M}}=\tau\left(1+\theta\left((1-\tau
\theta)I_t^S+I_t^M\right)\right)-1+\dfrac{1-\beta_3
P_{t+1}^R-\tau-\tau\widehat{d}
}{2}+\left(\beta_2-\beta_3/2\right)P_{t+1}^M,
\end{split}
\end{equation*}
and we get

\begin{equation}\label{a6}
\begin{split}
I_{t}^M=\dfrac{-\tau \theta I_t^S+\beta_3P_{t+1}^R
+(1+\widehat{d}-\tau)}{2 \tau
\theta}-\dfrac{(\beta_2-\beta_3/2)P_{t+1}^M}{\tau \theta}.
\end{split}
\end{equation}
Other constraints are
\begin{equation}\label{a8}
\begin{split}
P_{t}^M&=\dfrac{\partial H_t^{M}}{\partial
x_{t}}=2\widehat{{\delta}}x_t+\alpha
P_{t+1}^M+2\widehat{\widehat{\delta}} u_t,
\end{split}
\end{equation}

\begin{equation}\label{a9}
\begin{split}
u_{t+1}&=\dfrac{\partial H_t^{M}}{\partial P_{t+1}^R}=-\beta_3/2
I_t^M-\dfrac{\widehat{d}\beta_3}{2\tau \theta }-\dfrac{P_{t+1}^M
\beta_3^2}{2\tau \theta}+\alpha u_t.
\end{split}
\end{equation}

The equation of $I_t^M$ depends on $I_t^S$, and we obtain the
final equilibrium in the next section, at level two.
%-------------------------------------------------------------------------------------%
\label{sec:Mathematical Model: Level Two}\subsection{Mathematical
Model: Level Two}

At the previous level, the manufacturer's optimal function was
calculated by using a reaction function of the retailer. At this
level, the reaction functions of two followers (retailer and
manufacturer) are inserted into the objective function of the
leader (supplier), and we can find the final equilibrium point.\\
The problem facing the supplier is simply given by

\begin{equation*}
\begin{split}
&\arg \max_{I_t^S}
\sum_{t=1}^TP_t^Sq_t-cq_t+B_t^S(x_t)+T_t^S(I_t^S,I_t)-I_t^S+dI_t^{M},\\
\end{split}
\end{equation*}
subject to $x_{t+1}=\alpha x_t+\beta_1I_t^S+\beta_2I_t^{M}+\beta_3I_t^{R}$.\\

The Hamiltonian function of the supplier is defined by

\begin{equation}\label{a10}
\begin{split}
H_t^{S}&=J_t^{S}+P_{t+1}^{S}(x_{t+1})+\mu_t(P_t^M)+ u_t(p_t^R).
\end{split}
\end{equation}

Substitute (\ref{a0}) and (\ref{a6}) into (\ref{a10}), we get the
value of $I_t^S$, $x_{t+1}$ and $\mu_{t+1}$

\begin{equation}
\begin{split}
\dfrac{\partial H_t^{S}}{\partial I_{t}^S}=0,
\end{split}
\end{equation}

therefore
\begin{equation}
\begin{split}
I_t^S=\dfrac{(\beta_3/2+\beta_2-2\beta_1)}{\tau
\theta}P_{t+1}^S+\dfrac{(\beta_2-\beta_3/2)}{\tau
\theta}P_{t+1}^M+\dfrac{\beta_3}{2\tau\theta}P_{t+1}^R+\dfrac{3-3\tau
-\widehat{d}+2d}{2\tau \theta}.
\end{split}
\end{equation}

We have

\begin{equation}
\begin{split}
x_{t+1}&=\frac{\partial H_t^{S}}{\partial P_{t+1}^S},
\end{split}
\end{equation}
therefore we obtain
\begin{equation}
\begin{split}
x_{t+1}&=\alpha
x_t+(\beta_1-\beta_2/2-\beta_3/4)I_t^S+\dfrac{(-2\beta_2+\beta_3)(\beta_2-\beta_3/2)}{2\tau
\theta}P_{t+1}^M+\dfrac{(2\beta_2\beta3)-(3\beta_3^2)}{4\tau\theta}P_{t+1}^R\\
&+\dfrac{(2\beta_2-\beta_3)(1-\tau+\widehat{d})+2\beta_3(1-\tau)}{4\tau\theta}.
\end{split}
\end{equation}
We also have

\begin{equation}
\begin{split}
\mu_{t+1}&=\dfrac{\partial H_t^{S}}{\partial
P_{t+1}^M}=\alpha\mu_t-\frac{(\beta_2-\beta_3/2)}{2}I_t^S+\dfrac{(\beta_2-\beta_3/2)(\beta_3-2\beta_2)}{2\tau\theta}P_{t+1}^S\\
&-\frac{(\beta_2-\beta_3/2)d}{\tau\theta}.
\end{split}
\end{equation}
And we obtain

\begin{equation}
\begin{split}
P_{t}^S=\dfrac{\partial H_t^{S}}{\partial x_{t}}=2{\delta
x_t}+\alpha P_{t+1}^S+2{\widehat{\delta}}{u}_t.
\end{split}
\end{equation}

Since we use open-loop information, the structure variables depend
on the time variable and the initial state variables. The $x_1$ is
given initial parameter, $u_1 = 0$ and $\mu_1 = 0$. Furthermore,
the boundary condition are $P_{t+1}^R$= 0, $P_{t+1}^M$=0 and
$P_{t+1}^S$= 0.
%-----------------------------------------------------------------------------------%
\label{sec:Augmented Discrete Hamiltonian
Matrix}\subsection{Augmented Discrete Hamiltonian Matrix} In this
section for solving the optimal control problem formulated in
Section 3.1 and 3.2, we chose an algorithm given by Medanic and
Radojevic which is based on an augmented discrete Hamiltonian
matrix \cite{medanic}. First, we assume
\begin{equation*}
    \left[ {\begin{array}{c}
   \widetilde{x}_{t+1}\\
   \widetilde{P}_t\\
  \end{array} } \right]=\left[ {\begin{array}{cccc}
   A&B\\
   C&D\\
  \end{array} } \right]\left[ {\begin{array}{c}
   \widetilde{x}_{t}\\
   \widetilde{P}_{t+1}\\
  \end{array} } \right]+\left[ {\begin{array}{c}
   D\\
   E\\
  \end{array} } \right]=\left[ {\begin{array}{c}
   A \widetilde{x}_t+B \widetilde{P}_{t+1}+D\\
   C \widetilde{x}_t+A \widetilde{P}_{t+1}+E\\
  \end{array} } \right],
\end{equation*}
where $
  \widetilde{x}_{t+1}=
  \left[ {\begin{array}{c}
   x_{t+1}\\
   u_{t+1}\\
   \mu_{t+1}\\
  \end{array} } \right]
$ \quad and\quad $
  \widetilde{P}_{t+1}=
  \left[ {\begin{array}{c}
   p_{t+1}^S\\
   p_{t+1}^M\\
   p_{t+1}^R\\
  \end{array} } \right]$,

$A, B,$ and $C$ are $3\times3$ matrices, and $D$ and $E$ are
$3\times1$ matrices, such that
\begin{equation}
\begin{split}
  \widetilde{x}_{t+1}&=
  \left[ {\begin{array}{c}
   {x}_{t+1}\\
   {\mu}_{t+1}\\
   {u}_{t+1}\\
  \end{array} } \right]=\left[ {\begin{array}{c}
   A\widetilde{x}_t+B\widetilde{P}_{t+1}+D\\
  \end{array} } \right]=\left[ {\begin{array}{ccccccccc}
   a_{11}&a_{12}&a_{13}\\
   a_{21}&a_{22}&a_{23}\\
   a_{31}&a_{32}&a_{33}\\
  \end{array} } \right]\left[ {\begin{array}{c}
   x_t\\
   \mu_t\\
   u_t\\
  \end{array} } \right]+\left[ {\begin{array}{ccccccccc}
   b_{11}&b_{12}&b_{13}\\
   b_{21}&b_{22}&b_{13}\\
   b_{31}&b_{32}&b_{33}\\
  \end{array} } \right]\left[ {\begin{array}{c}
   P_{t+1}^S\\
   P_{t+1}^M\\
   P_{t+1}^R\\
  \end{array} } \right]\\
  &+\left[ {\begin{array}{c}
   d_1\\
   d_2\\
   d_3\\
  \end{array} } \right]=\left[ {\begin{array}{c}
   a_{11}x_t+a_{12}\mu_t+a_{13}u_t+b_{11}P^{S}_{t+1}+b_{12}P_{t+1}^{M}+b_{13}P_{t+1}^{S}+d_1\\
  a_{21}x_t+a_{22}\mu_t+a_{23}u_t+b_{21}P^{S}_{t+1}+b_{22}P_{t+1}^{M}+b_{23}P_{t+1}^{S}+d_2\\
  a_{31}x_t+a_{32}\mu_t+a_{33}u_t+b_{31}P^{S}_{t+1}+b_{32}P_{t+1}^{M}+b_{33}P_{t+1}^{S}+d_3\\
  \end{array} } \right].
  \end{split}
\end{equation}

 The boundary conditions are
 $
  \widetilde{x}_{1}=
  \left[ {\begin{array}{c}
   1\\
   0\\
   0\\
  \end{array} } \right]
$ \quad and\quad $
  \widetilde{P}_{T+1}=
  \left[ {\begin{array}{c}
   0\\
   0\\
   0\\
  \end{array} } \right]$.\\
$ A=\left[ {\begin{array}{ccccccccc}
   a_{11}&a_{12}&a_{13}\\
   a_{21}&a_{22}&a_{23}\\
   a_{31}&a_{32}&a_{33}\\
  \end{array} } \right],\\
$ where\\
 $a_{11}=\alpha,\quad a_{12}=0,\quad a_{13}=0,\quad
a_{21}=0,\quad a_{22}=\alpha,\quad a_{23}=0,\quad a_{31}=0,\quad
a_{32}=0,\quad
a_{33}=\alpha$.\\

$ B=\left[ {\begin{array}{ccccccccc}
   b_{11}&b_{12}&b_{13}\\
   b_{21}&b_{22}&b_{23}\\
   b_{31}&b_{32}&b_{33}\\
  \end{array} } \right]$,\\
  where

\begin{equation*}
b_{11}=\dfrac{(\beta_1-\beta_2/2-\beta_3/4)(-2\beta_1+\beta_2+\beta_3/2)}{\tau\theta}.
\end{equation*}

\begin{equation*}
b_{12}=\dfrac{(\beta_1-\beta_2/2-\beta_3/4)(\beta_2-\beta_3/2)}{\tau\theta}+\dfrac{(\beta_2-\beta_3/2)(-2\beta_2+\beta_3)}{2\tau\theta}.
\end{equation*}

\begin{equation*}
b_{13}=\dfrac{(\beta_1-\beta_2/2-\beta_3/4)(\beta_3)}{2\tau\theta}+\dfrac{(2\beta_2\beta_3)(-3\beta_3^2)}{4\tau\theta}.
\end{equation*}

\begin{equation*}
b_{21}=\dfrac{(\beta_2-\beta_3/2)(2\beta_1-3\beta_2+\beta_3/2)}{2\tau\theta}.
\end{equation*}
\begin{equation*}
b_{22}=-\frac{(\beta_2-\beta_3/2)^2}{2\tau\theta}.
\end{equation*}
\begin{equation*}
b_{23}=-\frac{\beta_3(\beta_2-\beta_3/2)}{4\tau\theta}.
\end{equation*}

\begin{equation*}
b_{31}=-\frac{\beta_3(-2\beta_1+\beta_2+\beta_3/2)}{4\tau\theta}.
\end{equation*}
\begin{equation*}
b_{32}=\frac{-\beta_3^2+3/2\beta_3(\beta_2-\beta_3/2)}{2\tau\theta}.
\end{equation*}
\begin{equation*}
b_{33}=\frac{-\beta_3^2}{8\tau\theta}.
\end{equation*}

$ D=\left[ {\begin{array}{ccc}
   d_{1}\\
   d_{2}\\
   d_{3}\\
  \end{array} } \right],$
  where

\begin{equation*}
\begin{split}
d_{1}&=\dfrac{(\beta_1-\beta_2/2-\beta_3/4)(3-3\tau-\widehat{d}+2d)}{2\tau\theta}
+\dfrac{(2\beta_2-\beta_3)(1-\tau+\widehat{d})+2\beta_3(1-\tau)}{4\tau\theta}.
\end{split}
\end{equation*}

\begin{equation*}
\begin{split}
d_{2}&=\dfrac{(-\beta_2-\beta_3/2)(6d-\widehat{d}-3\tau+3)}{4\tau\theta}.
\end{split}
\end{equation*}

\begin{equation*}
\begin{split}
d_{3}&=\dfrac{-\beta_3(-7\widehat{d}+2d-\tau+1)}{8\tau\theta}.
\end{split}
\end{equation*}

Similarly, we can get the values of the matrices $C$ and $E$
\begin{equation}
\begin{split}
  \widetilde{P}_{t}&=
  \left[ {\begin{array}{c}
   P_t^S\\
   P_t^M\\
   P_t^R\\
  \end{array} } \right]=\left[ {\begin{array}{ccccccccc}
   c_{11}&c_{12}&c{13}\\
   c_{21}&c_{22}&c{23}\\
   c_{31}&c_{32}&c{33}\\
  \end{array} } \right]\left[ {\begin{array}{c}
   x_t\\
  \mu_t\\
   u_t\\
  \end{array} } \right]+\left[ {\begin{array}{ccccccccc}
   a_{11}&a_{12}&a_{13}\\
   a_{21}&a_{22}&a_{23}\\
   a_{31}&a_{32}&a_{33}\\
  \end{array} } \right]\left[ {\begin{array}{c}
   P_{t+1}^S\\
   P_{t+1}^M\\
   P_{t+1}^R\\
  \end{array} } \right]+\left[ {\begin{array}{c}
   e_1\\
   e_2\\
   e_3\\
  \end{array} } \right]\\
  & =\left[ {\begin{array}{c}
  c_{11}x_t+c_{12}\mu_t+c_{13}u_t+a_{11}P_{t+1}^S+a_{12}P_{t+1}^{M}+a_{13}P_{t+1}^R+e_1\\
  c_{21}x_t+c_{22}\mu_t+c_{23}u_t+a_{21}P_{t+1}^S+a_{22}P_{t+1}^{M}+a_{23}P_{t+1}^R+e_2\\
  c_{31}x_t+c_{32}\mu_t+c_{33}u_t+a_{31}P_{t+1}^S+a_{32}P_{t+1}^{M}+a_{33}P_{t+1}^R+e_3\\
  \end{array} } \right].
  \end{split}
\end{equation}
Therefore
  $C=\left[
{\begin{array}{ccc}
   2\delta & 2\widehat{\delta}&2\widehat{\widehat{\delta}}\\
   2\widehat{\delta} &0&2\widehat{\widehat{\delta}}\\
   2\widehat{\widehat{\delta}}&0&0\\
  \end{array} } \right]$
and $ E=\left[ {\begin{array}{ccc}
   0\\
   0\\
   0\\
  \end{array} } \right].$

\label{sec:Resolution}\subsection{Resolution}

The above problem is solved by the sweep method \cite{Bryson}, by
assuming a linear relation between $\widetilde{p_t}$ and
$\widetilde{x_t}$
\\
\begin{equation}\label{319}
\widetilde{p}_k=S_k \widetilde{x}_k-g_k.
\end{equation}

Thus, we can compute

\begin{equation}\label{320}
\begin{split}
\widetilde{x_{k+1}}&=(I_{2*2}-BS_{k+1})^{-1}(A
\widetilde{x_t}-Bg_{k+1}+D).
\end{split}
\end{equation}
Then by substituting (\ref{319}) and (\ref{320}) into the
definition of $p_{k+1}$ as given by the augmented Hamiltonian
matrix, and equating both sides we finally get the difference
equations:

\begin{equation}
\begin{split}
S_{k}&=C+AS_{k+1}(I_{2*2}-BS_{k+1})^{-1}.
\end{split}
\end{equation}

\begin{equation}
\begin{split}
g_{k}&=AS_{k+1}(I_{2*2}-BS_{k+1})^{-1}+Bg_{k+1}-D + Ag_{k+1}-E.
\end{split}
\end{equation}

The boundary conditions are
$
  \widetilde{x}_{1}=
  \left[ {\begin{array}{c}
   x_1\\
   0\\
   0\\
  \end{array} } \right]
$ \quad and\quad $
  \widetilde{P}_{T+1}=
  \left[ {\begin{array}{c}
   0\\
   0\\
   0\\
  \end{array} } \right]$. And then $S_{T+1}=0_{3*3}$ and $g_{T+1}=0_{3*1}$.

From the boundary conditions we get $S_T=C$ and $g_T=E$. Once we
get the different values of $S_k$ and $g_k$ by the backward loop,
then the values of $\widetilde{x}_{t}$ and $\widetilde{p}_{t}$ are
computed by a forward loop. And, consequently we get the values of
$x_t$, $I_{t}^S$, $I_{t}^M$, $I_{t}^R$, $p_{t}^S$, $p_{t}^M$,
$p_{t}^R$, for all points in time.
%---------------------------------------------------------------------------------------%
\label{sec:Numerical Example}\section{Numerical Example} In this
section we provide a numerical example. We run the following
numerical simulations with mathematica 8. The results presented
here are obtained for the following values
of the parameters:\\
$a=6$, $w=3.8$, $c=2.4$, $q_t=100000$, $d=0.6$, $c=0.00001$,
$d=0.4$, $\widehat{d}=0.4$, $z=6$, $\theta=0.01$, and
 $\tau=0.2$.
We set $\beta_1=0.3$, $\beta_2=0.5$ and $\beta_3=0.8$; and
$\beta_1=0.3$, $\beta_2=0.5$. $B_t(x_t)=\delta x_t$, the potential
benefits firms obtain from
 taking social responsibility, such as increased demand, better
 reputation and so on. We set $\delta=0.2$, $\widehat{\delta}=0.2$ and
 $\widehat{\widehat{\delta}}=0.2$.
We assume that the time horizon is $T$=10. The initial level of
social responsibility is supposed to be $x_1$=1. We draw the
results of the equilibrium from our model, a three-stage Stackelberg
dynamic game.

\begin{figure}[ht!]
\begin{minipage}[b]{0.5\linewidth}
\centering
\includegraphics[width=\textwidth]{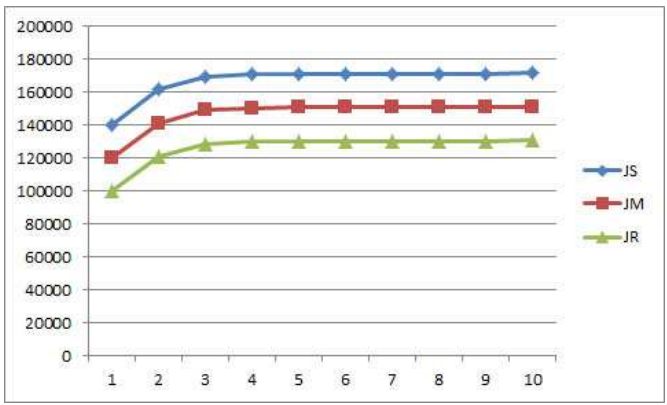}
\caption{Profits of supplier, manufacturer and retailer.}
\label{fig:figure3}
\end{minipage}
\end{figure}

\begin{figure}[ht!]
\begin{minipage}[b]{0.5\linewidth}
\centering
\includegraphics[width=\textwidth]{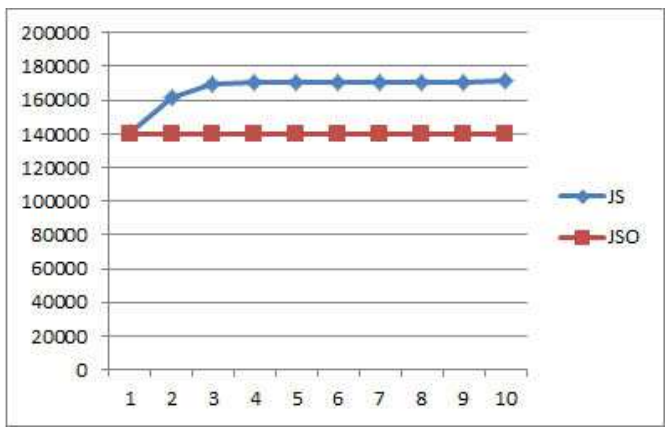}
\caption{Comparison of the supplier's profit, playing game one and
without playing any game.} \label{fig:figure4}
\end{minipage}
\end{figure}

\begin{figure}[ht!]
\begin{minipage}[b]{0.5\linewidth}
\centering
\includegraphics[width=\textwidth]{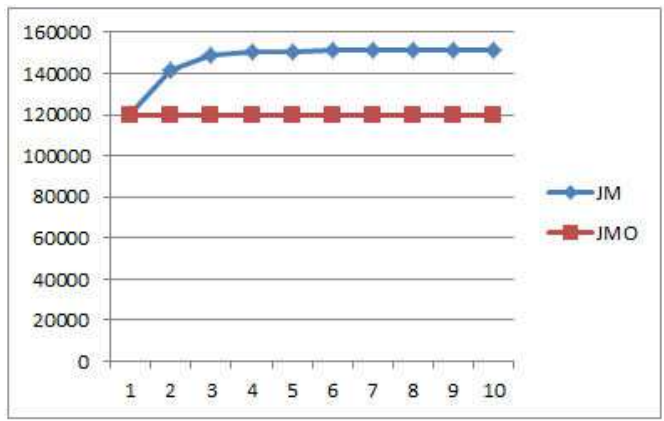}
\caption{Comparison of the manufacturer's profit, playing game one
and without playing any game.} \label{fig:figure5}
\end{minipage}
\end{figure}

\begin{figure}[!ht]
\begin{minipage}[b]{0.5\linewidth}
\centering
\includegraphics[width=\textwidth]{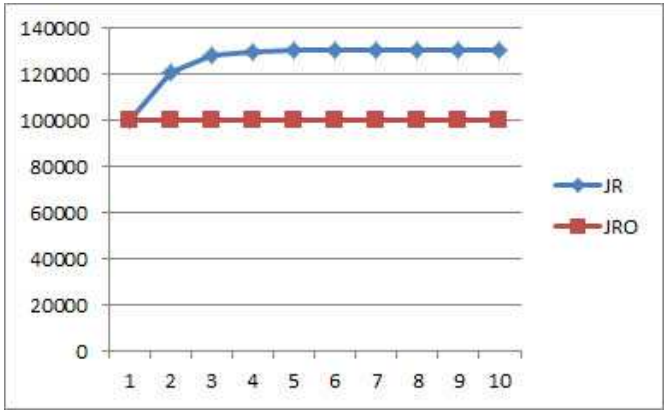}
\caption{Comparison of the retailer's profit, playing game one and
without playing any game.} \label{fig:figure6}
\end{minipage}
\end{figure}

\begin{figure}[!ht]
\begin{minipage}[b]{0.5\linewidth}
\centering
\includegraphics[width=\textwidth]{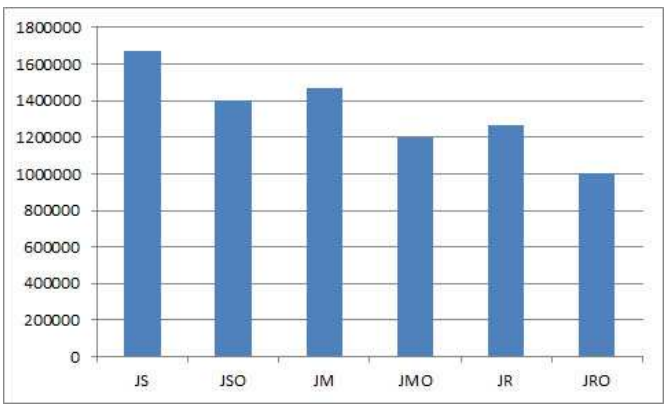}
\caption{Comparison of the cumulated profits of the member's of
supply chain, playing game one and without playing any game.}
\label{fig:figure7}
\end{minipage}
\end{figure}

The figure \ref{fig:figure3} shows the trend of profits from
periods one to ten in a Stackelberg game. $JS$ is the supplier's
profit, $JM$ is manufacturer's profit and $JR$ is retailer's
profit. We compare the profits of the supplier, manufacturer and
retailer over a time horizon, first while playing the game and
then, without playing the game. Figure \ref{fig:figure4} shows the
difference in supplier's profits when playing the game and without
playing. $JSO$ is supplier's profit without playing the game; $JS$
is supplier's profit when playing the game. As in the first graph,
the second and third one (figure \ref{fig:figure5},
\ref{fig:figure6}) shows the difference in manufacturer's profit
and retailer's profits when playing the game and without playing.
$JMO$ is manufacturer's profit without playing the game; $JM$ is
manufacturer's profit when playing the game.
 $JRO$ is retailer's profit without playing the game; $JR$
is retailer's profit when playing the game. Obviously, all of
players gain extra profit from playing the games. Figure
\ref{fig:figure7} compares the cumulated profits of the member's
of supply chain, playing game one and without playing game.

In sum, the supplier, manufacturer and retailer are motivated to
play the game because their respective benefit increases and the
supplier as the leader in the game earns more benefit than the
followers. Of course, this result is obtained with a very specific
dynamic game model. Another one may give different results.

%-------------------------------------------------------------------------------------------------%
\newpage
\label{sec:Conclusion}\section{Conclusion} In this paper we
investigated a decentralized three-tire supply chain consisting of
supplier, manufacturer and retailer with the aim of allocating CSR
to members of the supply chain system over time. We considered
two-level Stackelberg game consisting of two followers and one
leader. The members of a supply chain play games with each other
to maximize their own profits; thus, the model used was a
long-term co-investment game model. The equilibrium point in a
time horizon was determined at where the profit of supply chain's
members was maximized and CSR was implemented among members of the
supply chain. We applied control theory and used an algorithm
(augmented discrete Hamiltonian matrix) to obtain an optimal
solution for the dynamic game model. We presented a numerical
example and we found that, the benefits of the player increased
when they played the game.\\
\\
\textbf{Acknowledgments.} The present research was supported by
the MEDAlics, Research Center at Universit\`{a} per Stranieri
Dante Alighieri, Reggio Calabria, Italy.

%---------------------------------------------------------------------------------------%


\begin{thebibliography}{99}

\bibitem{Batabyal} Batabyal, A.A.: Consistency and optimality in a dynamic game of
pollution control II: Monopoly, Environmental and Resource
Economics, 8, 315-330 (1996).

\bibitem{Bryson} Bryson, A.E., Ho, Y.C.: Applied Optimal Control: Optimization,
Estimation and Control, John Wiley and Sons, New York (1975).

\bibitem{Cachon} Cachon, G.P., Zipkin, P.H.: Competitive and cooperative
inventory policies in a two-stage supply chain, Management
Science, 45, 936-953 (1999).

\bibitem{Carroll1} Carroll, A.B.: The pyramid of corporate social responsibility:
Toward the moral management of organizational stakeholders,
Business Horizons, 34, 39-48 (1991).

\bibitem{Carroll} Carroll, A.B.: A three-dimensional conceptual model of
corporate performance, Academy of Management Review, 4, 497-505
(1979).

\bibitem{Carter1} Carter, C.R., Kale, R., Grimm, C.M.: Environmental purchasing and
firm performance: An empirical investigation, Logistics and
Transportation Review, 36, 219-228 (2000).

\bibitem{Carter} Carter, C.R., Jennings, M.M.: Social responsibility and
supply chain relationships, Transportation Research, 38, 37-52
(2002).

\bibitem{Cetindamar} Cetindamar, D., Husoy, K.: Corporate social responsibility practices and environmentally responsible behavior: The case of the
united nations global compact, Journal of Business Ethics, 76,
163-176 (2007).

\bibitem{Dahlsrud} Dahlsrud, A.: How corporate social responsibility is defined: An analysis of 37 definitions, Corporate Social
Responsibility and Environmental Management, 15, 1-13 (2008).

\bibitem{Feibel} Feibel, B.J.: Investment Performance Measurement,  John Wiley and Sons,  New York, (2003).

\bibitem{He} He, X., Gutierrez, G.J., Sethi, S.P.: A survey of Stackelberg differential game models in
supply and marketing channels, Journal of Systems Science and
Systems Engineering, 16, 385-413 (2007).

\bibitem{Guti} He, X., Gutierrez, G., Sethi, S.P.: A review of Stackelberg differential game
models in supply chain management, Service Systems and Service
Management, (International Conference), 1, 9-11 (2007).

\bibitem{Hennet}Hennet, J.C., Arda, Y.: Supply chain coordination: A game-theory approach, Engineering Applications of Artificial Intelligence,
21,  399-405 (2008).

\bibitem{Hervani} Hervani, A.A., Helms, M.M., Sarkis, J.: Performance measurement for green supply management, Benchmarking: An International Journal, 12,  330-353 (2005).

\bibitem{Joyner} Joyner, B.E., Payne, D.: Evolution and implementation: A study of values, business ethics and corporate social responsibility, Journal
of Business Ethics, 41, 297-311 (2002).

\bibitem{Kamien} Kamien, M.I., Schwartz, N.L: Dynamic Optimization: The Calculus of Variations
and Optimal Control in Economics and Management,  Elsevier,
Amsterdam, (1991).

\bibitem{Mankiw} Mankiw, N.G.: Principles of Microeconomics, Mason, Ohio: Thomson/South-Western, (2004).

\bibitem{medanic} Medanic, J., Radojevic, D.: Multilevel Stackelberg strategies in linear-quadratic systems, Journal of Optimization Theory and Applications,
24, 485-497 (1978).

\bibitem{Orlitzky} Orlitzky, M.: Does firm size confound the relationship
between corporate social performance and firm financial
performance?, Journal of Business Ethics, 33, 167-180 (2001).

\bibitem{Preston} Preston, L.E.: The corporate
social-financial performance relationship: A Typology and
Analysis, Business and Society, 36, 419-429 (1997).

\bibitem{Sethi} Sethi, S.P.: Dimensions of corporate social responsibility: An Analytical Framework, California Management Review, 17, 58-64 (1975).

\bibitem{Shi} Shi, H.: A Game theoretic approach in green supply chain management, M.S. Thesis, University of Windsor, (2011).

\bibitem{Svensson} Svensson, G.: Aspects of sustainable supply chain management (SSCM): conceptual framework and empirical example, Supply Chain Management,
12, 262-266 (2007).

\bibitem{Tian} Tian, Y., Govindan, K., Zhu, Q.: A system dynamics model based on evolutionary game theory for green supply chain management diffusion among Chinese
manufacturers,  Journal of Cleaner Production, 80, 96-105 (2014).

\end{thebibliography}
\end{document}